\documentclass[preprint,12pt]{elsarticle}
\makeatletter
\def\ps@pprintTitle{%
  \let\@oddhead\@empty
  \let\@evenhead\@empty
  \let\@oddfoot\@empty
  \let\@evenfoot\@oddfoot
}
\makeatother

\usepackage{amssymb}
\usepackage{enumitem}
\usepackage{amsmath}
\usepackage{amsthm}
\newtheorem{theorem}{Theorem}[section]
\newtheorem{lemma}[theorem]{Lemma}
\newtheorem{corollary}[theorem]{Corollary}
\newtheorem{proposition}[theorem]{Proposition}
\usepackage{enumitem}
\usepackage{upgreek}
\usepackage{mathrsfs}

\usepackage{hyperref}
\usepackage{wasysym}
\usepackage[thinc]{esdiff}
\newtheorem{remark}{Remark}

\journal{}

\begin{document}

\begin{frontmatter}

%% Title, authors and addresses

%% use the tnoteref command within \title for footnotes;
%% use the tnotetext command for theassociated footnote;
%% use the fnref command within \author or \address for footnotes;
%% use the fntext command for theassociated footnote;
%% use the corref command within \author for corresponding author footnotes;
%% use the cortext command for theassociated footnote;
%% use the ead command for the email address,
%% and the form \ead[url] for the home page:
%% \title{Title\tnoteref{label1}}
%% \tnotetext[label1]{}
%% \author{Name\corref{cor1}\fnref{label2}}
%% \ead{email address}
%% \ead[url]{home page}
%% \fntext[label2]{}
%% \cortext[cor1]{}
%% \affiliation{organization={},
%%             addressline={},
%%             city={},
%%             postcode={},
%%             state={},
%%             country={}}
%% \fntext[label3]{}

\title{Drazin and g-Drazin invertibility of combinations of three Banach algebra elements}

%% use optional labels to link authors explicitly to addresses:
%% \author[label1,label2]{}
%% \affiliation[label1]{organization={},
%%             addressline={},
%%             city={},
%%             postcode={},
%%             state={},
%%             country={}}
%%
%% \affiliation[label2]{organization={},
%%             addressline={},
%%             city={},
%%             postcode={},
%%             state={},
%%             country={}}

\author[inst1]{Rounak Biswas}

\affiliation[inst1]{organization={Department of Mathematical And Computational Sciences, National Institute of Technology Karnataka},
            city={Surathkal}, 
            postcode={575025},
            state={Karnataka},
            country={India}.}
\author[inst1]{Falguni Roy}
\begin{abstract}
    Consider a complex unital Banach algebra $\mathcal{A}.$ For $x_1,x_2,x_3\in\mathcal{A},$ in this paper, we establish that under certain assumptions on $x_1,x_2,x_3$, Drazin (resp. g-Drazin) invertibility of any three elements among $x_1,x_2,x_3$ and $x_1+x_2+x_3\text{ }(\text{or }x_1x_2+x_1x_3+x_2x_3)$ ensure the Drazin (resp. g-Drazin) invertibility of the remaining one. As a consequence for two idempotents $p,q\in\mathcal{A},$ this result indicates the equivalence between Drazin (resp. g-Drazin) invertibility of $$\lambda_1p+\gamma_1q-\lambda_1pq+\lambda_2\left(pqp-(pq)^2\right)+\cdots+\lambda_m\left((pq)^{m-1}p-(pq)^m\right)$$ and $$\lambda_1-\lambda_1pq+\lambda_2\left(pqp-(pq)^2\right)+\cdots+\lambda_m\left((pq)^{m-1}p-(pq)^m\right),$$ 
    where $\gamma_1,\lambda_i\in\mathbb{C}$ for $i=1,2,\cdots,m,$ with $\lambda_1\gamma_1\neq0;$ which extend the work of J. Math. Anal. Appl.478 (2) (2019) 1163–1171. Furthermore, for $x_1,x_2$, we establish that the Drazin (resp. g-Drazin) invertibility of any two elements among $x_1,x_2$ and $x_1+x_2$ indicates the Drazin (resp. g-Drazin) invertibility of the remaining one, provided that $x_1x_2=\alpha(x_1+x_2)$ for some $\alpha\in\mathbb{C}$. Additionally, if it exists, we furnish a new formula to represent the Drazin (resp. g-Drazin) inverse of any element among $x_1,x_2$ and $x_1+x_2$,  by using the other two elements and their Drazin (resp. g-Drazin) inverse.
\end{abstract}

%%Graphical abstract
%\begin{graphicalabstract}
%\includegraphics{grabs}
%\end{graphicalabstract}

%%Research highlights
%\begin{highlights}
%\item Research highlight 1
%\item Research highlight 2
%\end{highlights}

\begin{keyword}
%% keywords here, in the form: keyword \sep keyword
Drazin inverse \sep Generalized Drazin inverse \sep Additive properties \sep Idempotents \sep Banach algebra
%% PACS codes here, in the form: \PACS code \sep code
%\PACS 0000 \sep 1111
%% MSC codes here, in the form: \MSC code \sep code
%% or \MSC[2008] code \sep code (2000 is the default)
\MSC 15A09 \sep 32A65 \sep 17C27  \sep 47A10
\end{keyword}

\end{frontmatter}
Email: rounak.207ma005@nitk.edu.in; royfalguni@nitk.edu.in
%\linenumbers

%% main text
\section{Introduction}

In this paper, $\mathcal{A}$ will represent an unital complex Banach algebra with the unit $1$.  An element $a\in\mathcal{A}$, is said to be Drazin invertible \cite{drazin1958pseudo} with the Drazin inverse $b\in\mathcal{A}$ $(\text{denoted by } a^d)$ if the following holds
\begin{equation}\label{Drazin_definition}
    ab=ba,\text{ }bab=b\text{ and }\left(a(1-ab)\right)^n=0\text{ for some }n\in\mathbb{N}. 
\end{equation} The smallest $n$ that satisfy (\ref{Drazin_definition}) is defined to be the Drazin index of $a$, which denoted by $i(a)$ and the collection of all Drazin invertible elements in $\mathcal{A}$ will be  represented by $\mathcal{A}^d$ . The Drazin inverse is a fundamental concept widely used across various domains of matrix analysis. It serves as a valuable tool for solving singular linear differential equations and difference equations \cite{campbell1976applications,campbell2009generalized} and holds significant importance in areas such as Markov chains \cite{meyer1980condition}, cryptography \cite{hartwig1981applications} and other areas \cite{behera2020further,kyrchei2014determinantal}. Koliha \cite{koliha1996generalized} introduces a generalized version of Drazin inverse, known as generalized Drazin (g-Drazin) inverse by using quasinilpotent \cite{harte1991quasinilpotents} elements (elements with spectral radius $0$) in  (\ref{Drazin_definition}). An element $a\in\mathcal{A}$ is said to be g-Drazin invertible with the g-Drazin inverse $b\in\mathcal{A}$ (denoted by $a^D$) if the following holds $$ab=ba,\text{ }bab=b\text{ and }a(1-ab)\in\mathcal{A}^{\text{qnil}},$$ where $\mathcal{A}^{\text{qnil}}$ represent the collection of all quasinilpotent elements in $\mathcal{A}.$
The collection of all g-Drazin invertible elements in $\mathcal{A}$ will be denoted by $\mathcal{A}^D.$

One of the most interesting problems in Drazin inverse theory is to study the Drazin (resp. g-Drazin) invertibility of $x_1+x_2$ when $x_1,x_2\in\mathcal{A}^d\text{ }(\text{resp. }\mathcal{A}^D).$ In 1958, Drazin \cite{drazin1958pseudo} proved that if $x_1,x_2\in\mathcal{A}^d$ where $x_1x_2=x_2x_1=0$, then $x_1+x_2\in\mathcal{A}^d$ with $(x_1+x_2)^d=x_1^d+x_2^d;$ and the same result is proved for g-Drazin inverse by Koliha \cite{koliha1996generalized}. In the case of matrices, just under the assumption $x_1x_2=0$, Hartwig et al. \cite{hartwig2001some} derived the representation for $(x_1+x_2)^d$ using $x_1,x_2,x_1^d$ and $x_2^d.$ This result is further generalized for g-Drazin inverse in \cite{gonzalez2004new}, within a Banach algebra under a weaker assumption. Over the last two decades, many researchers have worked on this additive problem of Drazin (resp. g-Drazin) inverse under different assumptions, some of which can be found in \cite{castro2009drazin,deng2009note,cvetkovic2013some,cvetkovic2011drazin,vivsnjic2015additive}.  

In the context of additive problems, if the choices of two elements from $\mathcal{A}$ are restricted to idempotents only, then investigating the Drazin invertibility of their sum reveals numerous intriguing findings. For two idempotents $p,q\in\mathcal{A},$ Koliha et al. \cite{koliha2012generalized} establish that the Drazin (resp. g-Drazin) invertibility of $\lambda_1 p+\gamma_1 q$ is equivalent to the Drazin (resp. g-Drazin) invertibility of $1-pq,$ where $\lambda_1, \gamma_1\in\mathbb{C}\setminus\{0\}.$ The proof of this result depends on the following equation satisfied by the idempotents $p$ and $q$
\begin{align*}
    \left(\lambda-\lambda_1(1-p)\right)&\left(\lambda-(\lambda_1 p+\gamma_1 q)\right)\left(\lambda-\gamma_1(1-q)\right)\notag\\&=\lambda\left((\lambda_1-\lambda)(\gamma_1-\lambda)-\lambda_1\gamma_1 pq\right),
\end{align*}where $\lambda\in\mathbb{C}.$ Recently, Benabdi and Barraa \cite{barraa2019drazin} proved the equivalence between the Drazin (resp. g-Drazin) invertibility of $\lambda_1 p+\gamma_1 q-\lambda_2 pq$ and $1-pq$ by using the equation \begin{align*}
    \left(\lambda-\lambda_1(1-p)\right)&\left(\lambda-(\lambda_1 p+\gamma_1 q-\lambda_2pq)\right)\left(\lambda-\gamma_1(1-q)\right)\notag\\&=\lambda\left(\lambda^2-\lambda(\lambda_1+\gamma_1)+\lambda_1\gamma_1+(\lambda_2\lambda-\lambda_1\gamma_1)pq\right),
\end{align*} where $\lambda_1,\gamma_1\in\mathbb{C}\setminus\{0\}$ and $\lambda_2\in\mathbb{C}.$ Motivated by these findings, in this paper, we present new additive results regarding the Drazin (resp. g-Drazin) invertibility, which are determined by the behaviour of the product \begin{equation}\label{product}
    (\lambda-x_1)(\lambda-x_2)(\lambda-x_3)
\end{equation} under some assumptions on $x_1,x_2,x_3$ where $x_1,x_2,x_3\in\mathcal{A}$ and $\lambda\in\mathbb{C}.$

Section \ref{sec_2} is dedicated to the study of the product (\ref{product}) while considering certain assumptions about $x_1,x_2,x_3$. As a result here we establish that if any three elements out of $x_1,x_2,x_3$ and $x_1+x_2+x_3\text{ }(\text{or }x_1x_2+x_1x_3+x_2x_3)$ are Drazin (resp. g-Drazin) invertible then the remaining one is also Drazin (resp. g-Drazin) invertible. Next, in Section \ref{idempotent}, these results are utilized to prove new additive results regarding Drazin (resp. g-Drazin) invertibility of combinations of two idempotents. In Section \ref{sec_3}, we investigate the Drazin (resp. g-Drazin) invertibility of $x_1,x_2$ and $x_1+x_2$ by considering the product $(\lambda-x_1)(\lambda-x_2).$ Moreover, here we establish a correlation between the representation of the Drazin (resp. g-Drazin) inverse of $x_1,x_2$ and $x_1+x_2$ whenever they exist. This recovers the scenario $x_1x_2=0$ examined by Hartwig et al. \cite{hartwig2001some}, and also confirms the converse aspect of their result.
\section{Additive results for three elements}\label{sec_2}
The aim of this section is to develop a method using the product (\ref{product}), which is useful for investigating the Drazin (resp. g-Drazin) invertibility of $x_1,x_2,x_3$, $x_1+x_2+x_3$ and $x_1x_2+x_1x_3+x_2x_3$ where $x_1,x_2,x_3\in\mathcal{A}.$ From now onwards, $\mathcal{A}^{\text{inv}}$ will denote the set of all invertible elements within $\mathcal{A}$, while sp$(a)$ and $\sigma(a)$ will refer to the spectral radius and spectrum of $a\in\mathcal{A},$ respectively. We begin with the following two lemmas for further use. 

\begin{lemma}\cite{shi2013drazin}\label{Drazin}
    Suppose $x\in\mathcal{A}.$ Then $x$ is Drazin invertible if and only if $0$ is not an essential singularity for the resolvent $(\lambda -x)^{-1}$ of $x$; and if $x\in\mathcal{A}^d$ then on the disc $\{\lambda:0<|\lambda|<(\text{sp}(x^d))^{-1}$\}, $(\lambda-x)^{-1}$ have the following Laurent series expansion,
    $$(\lambda-x)^{-1}=\sum_{n=1}^{i(x)}\frac{x^{n-1}}{\lambda^n}(1-xx^d)-\sum_{n=0}^{\infty}\lambda^n(x^d)^{n+1}.$$
\end{lemma}
\begin{lemma}\cite{koliha1996generalized}\label{Drazin_l1}
    Let $x\in\mathcal{A},$ then $x$ is g-Drazin invertible if and only if $0$ is not an accumulation point of $\sigma(x)$; and if $x$ is g-Drazin invertible then on the disc $\{\lambda:0<|\lambda|<(\text{sp}(x^D))^{-1}$\}, $(\lambda-x)^{-1}$ have the following Laurent series expansion,
    $$(\lambda-x)^{-1}=\sum_{n=1}^{\infty}\frac{x^{n-1}}{\lambda^n}(1-xx^D)-\sum_{n=0}^{\infty}\lambda^n(x^D)^{n+1}.$$
\end{lemma}

Now, we are ready to establish one of the main theorems of this section.

\begin{theorem}\label{3-Drazin}
    Suppose $x_1,x_2,x_3\in\mathcal{A}$ such that 
    \begin{align*}
        x_1x_2+x_1x_3+x_2x_3&=\alpha_1+\alpha_2(x_1+x_2+x_3)\text{ and }\\
        x_1x_2x_3&=\beta(x_1+x_2+x_3),
    \end{align*}       where $(\alpha_1,\alpha_2,\beta)\in\mathbb{C}^3\setminus\mathcal{M},$ with $\mathcal{M}=\{(\alpha_1,\alpha_2,0):\alpha_1,\alpha_2\in\mathbb{C}\text{ and }\alpha_1\alpha_2\neq0\}\cup \{(\alpha_1,0,0):\alpha_1\in\mathbb{C}\setminus\{0\}\}.$ Let $X$ denote the set $\{x_1,x_2,x_3,x_1+x_2+x_3\},$ then the Drazin invertibility of any three elements from the set $X$ indicates the Drazin invertibility of the remaining one.
     \begin{proof}
         Let $\lambda\in\mathbb{C},$ then 
      \begin{align}
          (\lambda-x_1)(\lambda-x_2)(\lambda-x_3)=&\text{ }\lambda^3-\lambda^2(x_1+x_2+x_3)+\lambda(x_1x_2+x_1x_3+x_2x_3)\notag\\&\text{ }-x_1x_2x_3\label{for_new}\\=&\text{ }\lambda^3+\lambda\alpha_1-(x_1+x_2+x_3)(\lambda^2-\lambda\alpha_2+\beta)\label{fromula}.
      \end{align}
         Now there are four ways to choose
         three elements from the set $X$. First, let us consider the case when $x_1,x_2,x_3$ are Drazin invertible. Then by Lemma \ref{Drazin}, for $0<|\lambda|<\text{ min }\{(\text{sp}(x_1^d))^{-1},(\text{sp}(x_2^d))^{-1},(\text{sp}(x_3^d))^{-1}\}$, $(\lambda-x_1), (\lambda-x_2),(\lambda -x_3)$ all are inveritible and $0$ is not an essential singularity for any of the resolvents $(\lambda-x_1)^{-1},(\lambda-x_2)^{-1},(\lambda-x_3)^{-1}.$ Furthermore for sufficiently small $\lambda$, $\lambda^2-\lambda\alpha_2+\beta$ is also invertible. Hence for sufficiently small $\lambda(\neq0)$, from (\ref{fromula}) we obtain 
         \begin{align}\label{inverseequation}
             \left(\frac{\lambda^3+\lambda\alpha_1}{\lambda^2-\lambda\alpha_2+\beta}-(x_1+x_2+x_3)\right)^{-1}=&\text{ }(\lambda^2-\lambda\alpha_2+\beta)(\lambda-x_3)^{-1}(\lambda-x_2)^{-1}\notag\\&\times(\lambda-x_1)^{-1}.
         \end{align} Suppose $\delta>0$ be small enough such as for $0<|\lambda|<\delta$, (\ref{inverseequation}) holds. Now let $D$ be a domain lying in $\mathbb{C}\setminus\{0\}$ which contain $\{\lambda:0<|\lambda|<\delta\}$ and let $f,g$ be two holomorphic functions on $D$ such that $(f(z^2))=z$, $(g(z^3))=z$ for $z\in D$. Next for sufficiently small $\mu$ substitute $$\lambda=-\frac{1}{3}\left(-\mu+\mathbf{A}+\frac{\Delta_0}{\mathbf{A}}\right)$$ in \ref{inverseequation}, where $\Delta_0=\mu^2-3\mu\alpha_2-3\alpha_1$, $\Delta_1=-2\mu^3+9\mu^2\alpha_2+\mu(9\alpha_1-27\beta)$ and $\mathbf{A}=g\left(\frac{\Delta_1+f(\Delta_1^2-4\Delta_0^3)}{2}\right)$. After this substitution, one can easily verify that $0$ is not an essential singularity for $(\mu-(x_1+x_2+x_3))^{-1},$ hence $x_1+x_2+x_3$ is Drazin invertible. Next, consider the case when $x_1,x_2,x_1+x_2+x_3$ from the set $X$ are Drazin invertible. Here also for small enough $\lambda(\neq 0)$, $(\lambda-x_1),(\lambda-x_2),(\lambda-(x_1+x_2+x_3))$ are invertible and $0$ is not an essential singularity for any of the resolvents $(\lambda-x_1)^{-1}$, $(\lambda-x_2)^{-1}$ and $(\lambda-(x_1+x_2+x_3))^{-1}.$ Since $(\alpha_1,\alpha_2,\beta)\in\mathbb{C}^3\setminus\mathcal{M},$  therefore for sufficiently small $\lambda(\neq0),$ from (\ref{fromula}) it follows
         \begin{align*}
             (\lambda-x_3)^{-1}=&\text{ }(\lambda^2-\lambda\alpha_2+\beta)^{-1}\left(\frac{\lambda^3+\lambda\alpha_1}{\lambda^2-\lambda\alpha_2+\beta}-(x_1+x_2+x_3)\right)^{-1}(\lambda-x_1)\\&\text{ }\times (\lambda-x_2).
         \end{align*} Thus $0$ is not an essential singularity for $(\lambda-x_3)^{-1}$, and hence $x_3$ is Drazin invertible. The other two possibilities can be verified in a similar way.
         
     \end{proof}
\end{theorem}
Similar to Theorem \ref{3-Drazin}, one can confirm the following outcome concerning g-Drazin invertibility through the utilization of Lemma \ref{Drazin_l1}.
\begin{theorem}\label{3-g-Drazin}
      Suppose $x_1,x_2,x_3\in\mathcal{A}$ such that 
    \begin{align*}
        x_1x_2+x_1x_3+x_2x_3&=\alpha_1+\alpha_2(x_1+x_2+x_3)\text{ and }\\
        x_1x_2x_3&=\beta(x_1+x_2+x_3),
    \end{align*}
    where $(\alpha_1,\alpha_2,\beta)\in\mathbb{C}^3\setminus\mathcal{M},$ with $\mathcal{M}=\{(\alpha_1,\alpha_2,0):\alpha_1,\alpha_2\in\mathbb{C}\text{ and }\alpha_1\alpha_2\neq0\}\cup \{(\alpha_1,0,0):\alpha_1\in\mathbb{C}\setminus\{0\}\}.$ Let $X$ denote the set $\{x_1,x_2,x_3,x_1+x_2+x_3\},$ then the g-Drazin invertibility of any three elements from the set $X$ indicates the g-Drazin invertibility of the remaining one.
\end{theorem}
Now, instead of satisfying the assumptions of Theorem \ref{3-Drazin} (resp. Theorem \ref{3-g-Drazin}) if $x_1,x_2,x_3$ satisfy the assumptions
\begin{align}\label{req_con}
       x_1+x_2+x_3&=\alpha_1+\alpha_2( x_1x_2+x_1x_3+x_2x_3)\text{ and }\notag\\
        x_1x_2x_3&=\beta(x_1x_2+x_1x_3+x_2x_3),
    \end{align}for $\alpha_1,\alpha_2,\beta\in\mathbb{C}$, then we have the following result:
    \begin{theorem}\label{3-gDTheoremTheoremTheorem-Drazin}
      Suppose $x_1,x_2,x_3\in\mathcal{A}$ satisfies the conditions $(\ref{req_con})$, and let $X$ represent the set $\{x_1,x_2,x_3,x_1x_2+x_1x_3+x_2x_3\}.$ Then the Drazin (resp. g-Drazin) invertibility of any three elements from the set $X$ indicates the Drazin (resp. g-Drazin) invertibility of the remaining one.
     \begin{proof}
         Since $x_1,x_2,x_3\in\mathcal{A},$ satisfies the assumptions (\ref{req_con}) here, then from (\ref{for_new}) we obtain
         \begin{align}\label{addition}
             (\lambda-x_1)(\lambda-x_2)(\lambda-x_3)=&\text{ }\lambda^3-\lambda^2\alpha_1\notag\\&\text{ }-(\lambda^2\alpha_2-\lambda+\beta)(x_1x_2+x_1x_3+x_2x_3).
         \end{align} Now if $x_1,x_2,x_3\in\mathcal{A}^d\text{ (resp. }\mathcal{A}^D)$ then for sufficiently small $\lambda(\neq0)$ from (\ref{addition}) we have 
         \begin{align}\label{neweq}
         \left(\frac{\lambda^3-\lambda^2\alpha_1}{\lambda^2\alpha_2-\lambda+\beta}-(x_1x_2+x_1x_3+x_2x_3)\right)^{-1}=& \text{ }(\lambda^2\alpha_2-\lambda+\beta)(\lambda-x_3)^{-1}\notag\\&\text{ }\times(\lambda-x_2)^{-1}(\lambda-x_1)^{-1}  .
         \end{align}
         Next, similar to the proof of Theorem \ref{3-Drazin}, using the functions $f\text{ and }g$, for small enough $\mu$ put the substitution $$\lambda=-\frac{1}{3}\left(-\alpha_1-\mu\alpha_2-\frac{1+\sqrt{-3}}{2}\mathbf{B}-\frac{2}{1+\sqrt{-3}}\frac{\Delta_{00}}{\mathbf{\textbf{B}}}\right),$$ in (\ref{neweq}), where $\Delta_{00}=\mu^2\alpha_2^2+\mu(2\alpha_1\alpha_2-3)+\alpha_1^2$, $\Delta_{10}=-2(\alpha_1+\mu\alpha_2)^3+\mu(9\alpha_1-27\beta+9\mu\alpha_2)$ and $\mathbf{B}=g\left(\frac{\Delta_{10}+f(\Delta_{10}^2-4\Delta_{00}^3)}{2}\right)$. From here onwards, the Drazin (resp. g-Drazin) invertibility of $x_1x_2+x_1x_3+x_2x_3$ follows similarly as in the proof of Theorem \ref{3-Drazin}. Employing the same method as demonstrated in the proof of Theorem \ref{3-Drazin}, one can verify the other three cases with the help of (\ref{addition}).
     \end{proof}
\end{theorem}
\begin{remark}\label{remrk}
 \normalfont For two idempotents $p,q\in\mathcal{A},$ some of the main results of \cite{koliha2012generalized} and \cite{barraa2019drazin} follows as a particular case of Theorem \ref{3-gDTheoremTheoremTheorem-Drazin}.
 \begin{enumerate}[label=(\roman*)]
     \item The equivalency among Drazin (resp. g-Drazin) invertibility of $\lambda_1p+\gamma_1q-\lambda_2pq$ and $1-pq$ established by Benabdi and Barraa  \cite[Theorem 2.2.]{barraa2019drazin}, follows by considering $x_1=\lambda_1(1-p),$ $x_2=\lambda_1p+\gamma_1q-\lambda_2pq$ and $x_3=\gamma_1(1-q)$ in Theorem \ref{3-gDTheoremTheoremTheorem-Drazin}, where $\lambda_1,\gamma_1\in\mathbb{C}\setminus\{0\}$ and $\lambda_2\in\mathbb{C}.$
     \item Similarly considering $x_1=\lambda_1(1-p),$ $x_2=\lambda_1p+\gamma_1q$ and $x_3=\gamma_1(1-q)$ in Theorem \ref{3-gDTheoremTheoremTheorem-Drazin} we obtain equivalency between the Drazin invertibility of $1-pq$ and $\lambda_1p+\gamma_1q$ for $\lambda_1,\gamma_1\in\mathbb{C}\setminus\{0\}$; which is one of the main result in \cite[Theorem 3.3.]{koliha2012generalized}.
 \end{enumerate}
\end{remark}
\begin{remark}
   \normalfont In Remark \ref{remrk}, using two idempotents we constructed $x_1,x_2$, and $x_3$, satisfying the assumptions of Theorem \ref{3-gDTheoremTheoremTheorem-Drazin}. This naturally raises the question: Can we generate further elements by employing two idempotents that meet the conditions of Theorem \ref{3-gDTheoremTheoremTheorem-Drazin}? Upon investigating this question, we find that $x_1=\lambda_1(1-p),x_2=\lambda_1p+\gamma_1q-\lambda_1pq+\lambda_2\left(pqp-(pq)^2\right)+\cdots+\lambda_m\left((pq)^{m-1}p-(pq)^m\right)$, and $x_3=\gamma_1(1-q)$ also satisfy the assumptions of Theorem \ref{3-gDTheoremTheoremTheorem-Drazin}. This selection of $x_1,x_2,x_3$ can be utilized to extend the work of \cite{koliha2012generalized} and \cite{barraa2019drazin}, which is elaborated in the next section.
\end{remark}
\section{Additivie results for idempotesnts}\label{idempotent}
One of the primary purposes of establishing new additive results mentioned in Section \ref{sec_2} is to utilize them to comprehend the  Drazin (resp. g-Drazin) invertibility when dealing with combinations of idempotents. The following result utilizes Theorem \ref{3-gDTheoremTheoremTheorem-Drazin} to demonstrate the equivalence between the Drazin (resp. g-Drazin) invertibility of two combinations of idempotents.
\begin{proposition}\label{prop}
    Let $p,q\in\mathcal{A}$ be two idempotents, and $\gamma_1,\lambda_i\in\mathbb{C}$ for $i=1,2,\cdots,m,$ with $\lambda_1\gamma_1\neq0.$ Then for $m\in\mathbb{N}$ $$\lambda_1p+\gamma_1q-\lambda_1pq+\lambda_2\left(pqp-(pq)^2\right)+\cdots+\lambda_m\left((pq)^{m-1}p-(pq)^m\right)$$ is Drazin (resp. g-Drazin) invertible if and only if $$\lambda_1-\lambda_1pq+\lambda_2\left(pqp-(pq)^2\right)+\cdots+\lambda_m\left((pq)^{m-1}p-(pq)^m\right),$$ is Drazin (resp. g-Drazin) invertible.
    \begin{proof}
        Choose \begin{align*}
            x_1=&\lambda_1(1-p),\\
            x_2=&\lambda_1p+\gamma_1q-\lambda_1pq+\lambda_2\left(pqp-(pq)^2\right)+\cdots+\lambda_m\left((pq)^{m-1}p-(pq)^m\right),\\
            x_3=&\gamma_1(1-q),
        \end{align*}  then $x_1,x_2$ and $x_3$ satisfies conditions (\ref{req_con}) with $\alpha_1=\gamma_1,\alpha_2=\frac{1}{\gamma_1}$ and $\beta=0.$ Hence the desired result follows from Theorem \ref{3-gDTheoremTheoremTheorem-Drazin}.
    \end{proof}   
\end{proposition}
\begin{corollary}
    Suppose $p,q\in\mathcal{A}$ are two idempotents, then for $m\in\mathbb{N}$ $$p+q-pq+pqp-(pq)^2+\cdots-(pq)^m$$ is Drazin (resp. g-Drazin) invertible if and only if $1-(pq)^m$ is Drazin (resp. g-Drazin) invertible.
    \begin{proof}
        As demonstrated by Koliha et al. \cite{koliha2012generalized}, corresponding to the idempotent  $p,$ $p$ and $q$ have the following matrix representation $$p=\begin{bmatrix}
            1 & 0\\
            0 & 0 \\
        \end{bmatrix}_p\text{ and }q=\begin{bmatrix}
            q_1 & q_2 \\
            q_3 & q_4 \\
        \end{bmatrix}_p.$$ Now, according to Proposition \ref{prop}, $p+q-pq+pqp-(pq)^2+\cdots-(pq)^m$ is Drazin (resp. g-Drazin) invertible if and only if \begin{align*}
            1-pq+pqp-\cdots-(pq)^m=\begin{bmatrix}
                1-q_1^m & -q_2-q_1q_1-\cdots-q_1^{m-1}q_2\\
                0 & 1 \\
            \end{bmatrix}_p
        \end{align*}is Drazin (resp. g-Drazin) invertible. From here, the required result can be verified using \cite[Lemma 2.4]{koliha2012generalized}.
    \end{proof}
\end{corollary}

Next, we recall the following lemma from \cite{koliha2012generalized}, which proves helpful in examining the Drazin (resp. g-Drazin) invertibility of a combination involving two idempotents. 
\begin{lemma}\cite{koliha2012generalized}\label{useful}
    Let $x\in\mathcal{A}$, and $f$ be a holomorphic function in an open neighbourhood of $\sigma(x),$ such that $f^{-1}(0)\cap \sigma(x)=\{\delta_1,\cdots,\delta_n\}$ for some $n\in\mathbb{N}$ and $\delta_k\in \mathbb{C}$ for $k=1,2\cdots n$. Then \begin{align*}
    f(x)\in \mathcal{A}^d&\iff \delta_k-x\in\mathcal{A}^d\text{ for all }k\in\{1,2,\cdots n\};\\
        f(x)\in \mathcal{A}^D&\iff \delta_k-x\in\mathcal{A}^D\text{ for all }k\in\{1,2,\cdots n\}.   
    \end{align*}
\end{lemma}
\begin{proposition}
    Suppose $p,q\in\mathcal{A},$ are two idempotents, then for $m\in\mathbb{N},$ the following are equivalent 
    \begin{enumerate}[label=\normalfont(\roman*)]
        \item\label{11} $p+q\in\mathcal{A}^d\text{ }(\text{resp. }A^D)$ and $p+q-1\in\mathcal{A}^d\text{ }(\text{resp. }A^D);$
        \item\label{22} $(pq)_m+(qp)_m\in\mathcal{A}^d\text{ }(\text{resp. }A^D)$;
        \item\label{33} $\begin{aligned}[t]
           \binom{m-1}{0}\left(pq+qp\right)&+\binom{m-1}{1}\left(pqp+qpq\right)+\cdots \\&+\binom{m-1}{m-1}\left((pq)_{m+1}+(qp)_{m+1}\right)\in\mathcal{A}^d\text{ }(\text{resp. }A^D);
        \end{aligned}$
        \item \label{44}$\begin{aligned}[t]
            \binom{m-1}{0}\left((pq)_{m+1}\right.&\left.+(qp)_{m+1}\right)+\binom{m-1}{1}\left((pq)_{m+2}+(qp)_{m+2}\right)+\cdots\\& +\binom{m-1}{m-1}\left((pq)_{2m}+(qp)_{2m}\right)\in\mathcal{A}^d\text{ }(\text{resp. }A^D);
        \end{aligned}$
    \end{enumerate}where $\binom{m}{k}=\frac{m!}{(k!)(m-k)!}$, $(pq)_k=\underbrace{pqpqp\cdots}_{k\text{ times}},$ and $(qp)_k=\underbrace{qpqpq\cdots}_{k\text{ times}}$, for $k\in\mathbb{N}.$
    \begin{proof}
      \ref{11}$\iff$\ref{22} Consider the holomorphic function $f(\alpha)=\alpha(\alpha-1)^{m-1}$, where $\alpha\in\mathbb{C}.$ Then a simple calculation confirms that $f(p+q)=(pq)_m+(qp)_m, $ from here, the required result follows by Lemma \ref{useful}.
    \\ Similarly the equivalency between \ref{11}, \ref{33} and \ref{11}, \ref{44} can be verified by considering the functions $f(\alpha)=\alpha^m(\alpha-1)$ and $f(\alpha)=\alpha^m(\alpha-1)^m,$ respectively.
    \end{proof}
\end{proposition}
\section{Additive results for two elements}\label{sec_3}
Instead of considering three factors in (\ref{product}), if we consider just two factors, then the simplicity of the product allows us to investigate the representation for Drazin inverse also.
In this section, we will study the product $(\lambda-x_1)(\lambda-x_2)$ while assuming that $x_1x_2=\alpha(x_1+x_2)$ where $x_1,x_2\in\mathcal{A}$ and $\alpha\in \mathbb{C}.$ 

\begin{theorem}\label{2deg_Drazin}
    Let $x_1,x_2\in\mathcal{A},$ such that $x_1x_2=\alpha (x_1+x_2)$, for some $\alpha \in \mathbb{C}.$
    \begin{enumerate}[label=\normalfont(\roman*)]
        \item If $x_1,x_2\in\mathcal{A}^d$ then $(x_1+x_2)\in\mathcal{A}^d$ and 
        \begin{align}\label{fromula1}
            (x_1+x_2)^d=&\sum_{n=1}^{i(x_2)}(1-x_2x_2^d)x_2^{n-1}(x_1^d)^n+\sum_{n=1}^{i(x_1)}(x_2^d)^nx_1^{n-1}(1-x_1x_1^d)\notag\\&-2\alpha\left(\sum_{n=1}^{i(x_2)}(1-x_2x_2^d)x_2^{n-1}(x_1^d)^{n+1}\right.\notag\\&\left.+\sum_{n=1}^{i(x_1)}(x_2^d)^{n+1}x_1^{n-1}(1-x_1x_1^d)\right).
        \end{align}
         \item If $x_1+x_2\in \mathcal{A}^d$ and $x_1\text{ }(\text{or }x_2)\in\mathcal{A}^d$ then $x_2\text{ }(\text{or }x_1)\in\mathcal{A}^d,$ and \begin{align}\label{converse_formula}
             x_2^d=\begin{cases}
                 (x_1+x_2)^d-\left((x_1+x_2)^d\right)^2x_1,\text{ }\text{  if }\alpha=0\\\frac{1}{\alpha}(x_1+x_2)^dx_1,\text{ }\text{ if }\alpha\neq 0          
             \end{cases}.
         \end{align}
    \end{enumerate}
    \begin{proof}
     \begin{enumerate}[label=(\roman*)]
         \item Suppose $x_1,x_2\in\mathcal{A}^d$ then for $\lambda \in \mathbb{C}$ we have $(\lambda -x_1),(\lambda -x_2)\in\mathcal{A}^{\text{inv}}$ whenever $0<|\lambda|<\text{ min }\{(\text{sp}(x_1^d))^{-1},(\text{sp}(x_2^d))^{-1}\},$ and 
         \begin{align}
             (\lambda-x_1)^{-1}&=\sum_{n=1}^{i(x_1)}\frac{x_1^{n-1}}{\lambda^n}(1-x_1x_1^d)-\sum_{n=0}^{\infty}\lambda^n(x_1^d)^{n+1}; \label{x_1res}\\
             (\lambda-x_2)^{-1}&=\sum_{n=1}^{i(x_2)}\frac{x_2^{n-1}}{\lambda^n}(1-x_2x_2^d)-\sum_{n=0}^{\infty}\lambda^n(x_2^d)^{n+1} ,\label{x_2res}
        \end{align}
         
      Now since $x_1x_2=\alpha (x_1+x_2)$, therefore for $\lambda \in \mathbb{C}\setminus \{\alpha\}$ we have 
         \begin{align}\label{2nd}
             (\lambda -x_1)(\lambda -x_2)&=(\lambda-\alpha)\left(\frac{\lambda^2}{\lambda-\alpha}-(x_1+x_2)\right). 
         \end{align} Next, for sufficiently small $\mu$ substituting $\lambda=\frac{\mu+\sqrt{\mu^2-4\mu \alpha}}{2}$ in (\ref{2nd}), we get
         \begin{align}\label{2nd_mu}
            \left(\frac{\mu+\sqrt{\mu^2-4\mu \alpha}}{2}-\alpha\right)^{-1}& \left(\frac{\mu+\sqrt{\mu^2-4\mu \alpha}}{2}-x_1\right)\notag\\&\times\left(\frac{\mu+\sqrt{\mu^2-4\mu \alpha}}{2}-x_2\right)\notag\\=&\left(\mu-(x_1+x_2)\right).
         \end{align}
    Here if $\mu(\neq0)$ is sufficiently small, then $$0<|\frac{\mu+\sqrt{\mu^2-4\mu \alpha}}{2}|<\text{ min }\{(\text{sp}(x_1^d))^{-1},(\text{sp}(x_2^d))^{-1}\},$$ thus by (\ref{x_1res}) and (\ref{x_2res}), for sufficiently small $\mu,$ $0$ is not an essential singularity of $\left(\frac{\mu+\sqrt{\mu^2-4\mu \alpha}}{2}-x_1\right)^{-1}$ and $\left(\frac{\mu+\sqrt{\mu^2-4\mu \alpha}}{2}-x_2\right)^{-1}$, respectively. Therefore from (\ref{2nd_mu}), it follows that for $\mu$ sufficiently small 
     \begin{align}\label{inveq1}
         (\mu-(x_1+x_2))^{-1}=&\left(\frac{\mu+\sqrt{\mu^2-4\mu \alpha}}{2}-\alpha\right)\left(\frac{\mu+\sqrt{\mu^2-4\mu \alpha}}{2}-x_2\right)^{-1}\notag\\&\times\left(\frac{\mu+\sqrt{\mu^2-4\mu \alpha}}{2}-x_1\right)^{-1};
     \end{align}and $0$ is not an essential singularity of $(\mu-(x_1+x_2))^{-1}$. Hence we obtain $x_1+x_2\in\mathcal{A}^d,$ moreover by \cite[Remark 2.1.]{shi2013drazin} we have 
     \begin{equation}\label{inteq1}
         (x_1+x_2)^d=-\frac{1}{2\pi i}\oint_{\Gamma}\frac{1}{\mu}\left(\mu-(x_1+x_2)\right)^{-1}d\mu;
     \end{equation}
      where $\Gamma=\{\mu\in\mathbb{C}:|\frac{\mu+\sqrt{\mu^2-4\mu \alpha}}{2}|=\epsilon\}$ and $0<\epsilon<\text{ min }\{(\text{sp}(x_1^d))^{-1},(\text{sp}(x_2^d))^{-1}\}.$ 
     Now using (\ref{inveq1}) it follows  \begin{align}
         \oint_{\Gamma}\frac{1}{\mu}(\mu-(x_1&+x_2)^{-1})d\mu\notag\\=&\oint_{\Gamma}\left(\frac{\mu+\sqrt{\mu^2-4\mu \alpha}-2\alpha}{2\mu}\right)\left(\frac{\mu+\sqrt{\mu^2-4\mu \alpha}}{2}-x_2\right)^{-1}\notag\\&\times\left(\frac{\mu+\sqrt{\mu^2-4\mu \alpha}}{2}-x_1\right)^{-1}d\mu.\label{subs}
     \end{align} Next substituting $\gamma=\frac{\mu+\sqrt{\mu^2-4\mu \alpha}}{2}$ in (\ref{subs}) we obtain
     \begin{align}\label{intrep}
         \oint_{\Gamma}\frac{1}{\mu}\left(\mu-(x_1+x_2)^{-1}\right)d\mu=&\oint_{|\gamma|=\epsilon}\left(1-\frac{2\alpha}{\gamma}\right)(\gamma-x_2)^{-1}(\gamma-x_1)^{-1}d\gamma\notag\\=&\text{ }2\pi i\left(-\sum_{n=1}^{i(x_2)}(1-x_2x_2^d)x_2^{n-1}(x_1^d)^n\right.\notag\\&\text{ }-\sum_{n=1}^{i(x_1)}(x_2^d)^nx_1^{n-1}(1-x_1x_1^d)\notag\\&\text{ }+2\alpha\left(\sum_{n=1}^{i(x_2)}(1-x_2x_2^d)x_2^{n-1}(x_1^d)^{n+1}\right.\notag\\&\left.\left.\text{ }+\sum_{n=1}^{i(x_1)}(x_2^d)^{n+1}x_1^{n-1}(1-x_1x_1^d)\right)\right). 
     \end{align} Hence, the required representation (\ref{fromula1}) of $(x_1+x_2)^d$ follows from (\ref{inteq1}) and (\ref{intrep}).
     \item Next, suppose $x_1,x_1+x_2\in\mathcal{A}^d,$ then for small enough $\lambda(\neq0)$; $\lambda-\alpha,$ $\lambda-x_1$ and $\frac{\lambda^2}{\lambda-\alpha}-(x_1+x_2)$ are invertible. Hence for sufficiently small $\lambda$ from (\ref{2nd}) we obtain 
     \begin{align}\label{idk}
         (\lambda-x_2)^{-1}=&\frac{1}{\lambda-\alpha}\left(\frac{\lambda^2}{\lambda-\alpha}-(x_1+x_2)\right)^{-1}(\lambda-x_1),
     \end{align} where $0$ is not an essential singularity for the resolvent of $x_1+x_2$ i.e., $\left(\frac{\lambda^2}{\lambda-\alpha}-(x_1+x_2)\right)^{-1}$. This indicates that $0$ is not an essential singularity for $(\lambda-x_2)^{-1}$ as well, hence $x_2$ is Drazin invertible. Now let us assume that $\alpha\neq0,$ then (\ref{idk}) gives
     \begin{align*}
         (\lambda-x_2)^{-1}=&-\frac{1}{\alpha}\left(\sum_{k=0}^\infty\left(\frac{\lambda}{\alpha}\right)^k\right)\left(\sum_{n=1}^{i(x_1+x_2)}\left(\frac{1}{\lambda}-\frac{\alpha}{\lambda^2}\right)^n(x_1+x_2)^{n-1}\right.\\&\times\left(1-(x_1+x_2)(x_1+x_2)^d\right)\\&\left.-\sum_{n=0}^\infty\left(\frac{\lambda^2}{\lambda-\alpha}\right)^n\left((x_1+x_2)^d\right)^{n+1}\right)(\lambda-x_1),
     \end{align*}using this expression of $(\lambda-x_2)^{-1}$ and \cite[Remark 2.1.]{shi2013drazin} one can validate that $$x_2^d=\frac{1}{\alpha}(x_1+x_2)^dx_1.$$ Next if $\alpha=0$ then (\ref{idk}) becomes 
     $$(\lambda-x_2)^{-1}=\frac{1}{\lambda}\left(\lambda-(x_1+x_2)\right)^{-1}(\lambda-x_1),$$from here one can easily verify that $x_2^d=(x_1+x_2)^d-\left((x_1+x_2)^d\right)^2x_1$ as desired.
    \end{enumerate} 
    \end{proof}
\end{theorem}
\begin{corollary}\label{coro}
     Let $x_1,x_2\in\mathcal{A},$ such that $x_1x_2=0$.
    \begin{enumerate}[label=\normalfont(\roman*)]
        \item\label{1k} If $x_1,x_2\in\mathcal{A}^d$ then $(x_1+x_2)\in\mathcal{A}^d$ and \begin{align}\label{fromula}
            (x_1+x_2)^d=&\sum_{n=1}^{i(x_2)}(1-x_2x_2^d)x_2^{n-1}(x_1^d)^n+\sum_{n=1}^{i(x_1)}(x_2^d)^nx_1^{n-1}(1-x_1x_1^d)\notag.
        \end{align}
         \item\label{2k} If $x_1+x_2\in \mathcal{A}^d$ and $x_1\text{ }(\text{or }x_2)\in\mathcal{A}^d$ then $x_2\text{ }(\text{or }x_1)\in\mathcal{A}^d,$ where $$x_2^d=(x_1+x_2)^d-\left((x_1+x_2)^d\right)^2x_1.$$
    \end{enumerate}
    \begin{proof}
       This is a direct consequence of Theorem \ref{2deg_Drazin} when $\alpha=0.$
    \end{proof}
\end{corollary}
\begin{remark}
    \normalfont Result \ref{1k} of Corollary \ref{coro} is due to  Hartwig et al. \cite{hartwig2001some}. However, their work did not address the complementary aspect of this result, i.e., if $x_1x_2=0$ and $x_1,x_1+x_2\in\mathcal{A}^d$, then what conclusions can be made regarding the Drazin invertibility of $x_2.$ Our result \ref{2k} of Corollary \ref{coro} addresses this scenario. Recent work of Zou \cite{zou2024some} also addresses this converse aspect of Hartwig et al. \cite{hartwig2001some} results in a different approach. But the method used in \cite{zou2024some} is unable to give any representation formula for $x_2^d$ using $x_1,$ $x_1+x_2$ and their Drazin inverse under the assumption $x_1x_2=0$.
\end{remark}
Similar to Theorem \ref{2deg_Drazin}, the following result can be established, which addresses g-Drazin invertibility within the same context.
\begin{theorem}
    Let $x_1,x_2\in\mathcal{A},$ such that $x_1x_2=\alpha (x_1+x_2)$, for some $\alpha \in \mathbb{C}.$
    \begin{enumerate}[label=\normalfont(\roman*)]
        \item If $x_1,x_2\in\mathcal{A}^D$ then $(x_1+x_2)\in\mathcal{A}^D$ and 
        \begin{align*}
            (x_1+x_2)^D=&\sum_{n=1}^{\infty}(1-x_2x_2^D)x_2^{n-1}(x_1^D)^n+\sum_{n=1}^{\infty}(x_2^D)^nx_1^{n-1}(1-x_1x_1^D)\notag\\&-2\alpha\left(\sum_{n=1}^{\infty}(1-x_2x_2^D)x_2^{n-1}(x_1^D)^{n+1}\right.\notag\\&\left.+\sum_{n=1}^{\infty}(x_2^D)^{n+1}x_1^{n-1}(1-x_1x_1^D)\right).
        \end{align*}
         \item If $x_1+x_2\in \mathcal{A}^D$ and $x_1\text{ }(\text{or }x_2)\in\mathcal{A}^D$ then $x_2\text{ }(\text{or }x_1)\in\mathcal{A}^D,$ and \begin{align*}
             x_2^D=\begin{cases}
                 (x_1+x_2)^D-\left((x_1+x_2)^D\right)^2x_1,\text{ }\text{ if }\alpha=0\\\frac{1}{\alpha}(x_1+x_2)^Dx_1,\text{ }\text{ if }\alpha\neq 0,          
             \end{cases}.
         \end{align*}
    \end{enumerate}
\end{theorem}

%% If you have bibdatabase file and want bibtex to generate the
%% bibitems, please use
%%
 %\bibliographystyle{elsarticle-num} 
 %\bibliography{cas-refs}

%% else use the following coding to input the bibitems directly in the
%% TeX file.

% \begin{thebibliography}{00}

% %% \bibitem{label}
% %% Text of bibliographic item

% \bibitem{}

% \end{thebibliography}
\end{document}